\documentclass[12pt,reqno]{amsart}

\usepackage{amssymb, amsmath, amsfonts, latexsym}
\usepackage{enumerate}
\usepackage{verbatim}
\usepackage{color}
\newtheorem{thm}{Theorem}[section]
\newtheorem{cor}[thm]{Corollary}
\newtheorem{lem}[thm]{Lemma}
\newtheorem{prop}[thm]{Proposition}

\newtheorem{remarks}[thm]{Remark}

\theoremstyle{definition}
\newtheorem{defn}{Definition}[section]

\numberwithin{equation}{section} \theoremstyle{remark}

\begin{document}

\title[Spectral gap and weighted Poincar\'e inequalities]{\bf A note on spectral gap and weighted Poincar\'e inequalities for some one-dimensional diffusions}

\author{Michel Bonnefont}
\address{Michel BONNEFONT \\ Institut de Math\'ematiques de Bordeaux, Universit\'e  de Bordeaux, 351 cours de la Lib\'eration,
33405 Talence, France.
}
\thanks{MB is partially supported by the French ANR-12-BS01-0013-02 HAB project}
\email{michel.bonnefont@math.u-bordeaux.fr}
\urladdr{http://www.math.u-bordeaux.fr/~mibonnef/}

\author{Ald\'eric Joulin}
\address{Ald\'eric JOULIN (corresponding author) \\ Universit\'e de Toulouse,
Institut National des Sciences Appliqu\'ees,
Institut de Math\'ematiques de Toulouse,
F-31077 Toulouse,
France.
}
\thanks{AJ is partially supported by the French ANR-2011-BS01-007-01 GEMECOD and ANR-12-BS01-0019 STAB projects}
\email{ajoulin@insa-toulouse.fr}
\urladdr{http://www.math.univ-toulouse.fr/~joulin/}

\author{Yutao Ma}
\address{Yutao MA\\ School of Mathematical Sciences $\&$ Lab. Math. Com. Sys., Beijing Normal University, 100875 Beijing, China.}
\thanks{YM is partially supported by NSFC 11101040, YETP0264, 985 Projects and the Fundamental Research Funds for the Central Universities.}
\email{mayt@bnu.edu.cn}
\urladdr{http://math.bnu.edu.cn/~mayt/}

\date{}

\newcommand{\aaa}{\mathbb{A}}
\newcommand{\bb}{\mathbb{B}}
\newcommand{\cc}{\mathbb{C}}
\newcommand{\dd}{\mathbb{D}}
\newcommand{\ee}{\mathbb{E}}
\newcommand{\mm}{\mathbb{M}}
\newcommand{\nn}{\mathbb{N}}
\newcommand{\rr}{\mathbb{R}}
\newcommand{\pp}{\mathbb{P}}
\newcommand{\qq}{\mathbb{Q}}
\newcommand{\ttt}{\mathbb{T}}
\newcommand{\zz}{\mathbb{Z}}

\def\Lip{{\mathrm{{\rm Lip}}}}
\def\Var{{\mathrm{{\rm Var}}}}
\def\AA{\mathcal A}
\def\BB{\mathcal B}
\def\CC{\mathcal C}
\def\DD{\mathcal D}
\def\FF{\mathcal F}
\def\EE{\mathcal E}
\def\II{\mathcal I}
\def\JJ{\mathcal J}
\def\KK{\mathcal K}
\def\LL{\mathcal L}
\def\MM{\mathcal M}
\def\NN{\mathcal N}
\def\SS{\mathcal S}
\def\TT{\mathcal T}
\def\WW{\mathcal W}

\def\vep{\varepsilon}
\def\<{\langle}
\def\>{\rangle}
\def\dsubset{\subset\subset}

\def\beq{\begin{equation}}
\def\nneq{\end{equation}}

\def\bdef{\begin{defn}}
\def\ndef{\end{defn}}

\def\bthm{\begin{thm}}
\def\nthm{\end{thm}}

\def\bprop{\begin{prop}}
\def\nprop{\end{prop}}

\def\brmk{\begin{remarks}}
\def\nrmk{\end{remarks}}

\def\bexa{\begin{exa}}
\def\nexa{\end{exa}}

\def\blem{\begin{lem}}
\def\nlem{\end{lem}}

\def\bcor{\begin{cor}}
\def\ncor{\end{cor}}

\def\bexe{\begin{exe}}
\def\nexe{\end{exe}}

\def\bprf{\begin{proof}}
\def\nprf{\end{proof}}

\def\lar{\overleftarrow}
\def\rar{\overrightarrow}

\def\bdes{\begin{description}}
\def\ndes{\end{description}}

\def\Cov{{\rm Cov}}

\keywords{Spectral gap; Diffusion operator; Weighted Poincar\'e inequality.}


\maketitle

\begin{abstract}
We present some classical and weighted Poincar\'e inequalities for some one-dimensional probability measures. This work is the one-dimensional counterpart of a recent study achieved by the authors for a class of spherically symmetric probability measures in dimension larger than 2. Our strategy is based on two main ingredients: on the one hand, the optimal constant in the desired weighted Poincar\'e inequality has to be rewritten as the spectral gap of a convenient Markovian diffusion operator, and on the other hand we use a recent result given by the two first authors, which allows to estimate precisely this spectral gap. In particular we are able to capture its exact value for some examples.
\end{abstract}

\section{Introduction}
Let $\MM$ be a smooth connected Riemannian manifold which is complete for the associated Riemannian metric $g$. Denote $\mu$ the probability measure on $\MM$ with density (with respect to the volume element on $\MM$) proportional to $e^{-V}$, where $V$ is some nice potential, and let $\sigma$ be some smooth function satisfying convenient assumptions. Let $\CC ^\infty (\MM)$ be the space of real-valued smooth functions on $\MM$ and let $\CC _0 ^\infty (\MM)$ be the subspace of $\CC ^\infty (\MM)$ consisting of compactly supported functions. We say that $\mu$ satisfies a weighted Poincar\'e inequality with constant $\lambda > 0$ and weight $\sigma^2$ if for any $f\in \CC _0 ^\infty (\MM)$, the following inequality holds:
\beq
\label{modipoin}
\lambda \Var _{\mu}(f)\le \int_\MM \sigma^2 \vert \nabla _\MM f \vert ^2 d\mu ,
\nneq
where $\nabla _\MM$ is the Riemannian gradient and $\vert \cdot \vert $ is the norm with respect to the metric $g$. Here $\Var _{\mu}(f)$ is the variance of $f$ with respect to $\mu$, that is
$$
\Var _{\mu}(f) := \int_\MM f^2 d\mu - \left( \int_\MM fd\mu \right) ^2 .
$$
Such a functional inequality reduces to the classical Poincar\'e inequality when there is no weight or, in other words, $\sigma \equiv 1$. Recently, the question whether a probability measure satisfies such a functional inequality has attracted a lot of attention, cf. for instance the papers \cite{BCG,miclo,BLe,CGGR,gozlan,bj,bjm}. Among the potential applications of these Poincar\'e type inequalities, one of the most relevant and interesting features is given in terms of measure concentration, large deviations and tail estimates of Lipschitz functionals. The interested reader is referred to the set of notes of Ledoux \cite{ledoux_berlin,Ledoux} for an introduction on the rich theory of functional inequalities and their basic (and less basic) properties, and also to the book of Chen \cite{chen_book} in which, among other situations, the one-dimensional case is carefully emphasized. \vspace{0.1cm}

An alternative probabilistic point of view on these functional inequalities consists in regarding the right-hand-side of \eqref{modipoin} as an energy term derived from a convenient (essentially self-adjoint) Markovian generator. More precisely, the measure $\mu$ is invariant for many reversible Markovian dynamics and among them, the diffusion operators we can consider are of the following form: for any $f\in \CC _0 ^\infty (\MM)$,
$$
\LL_{\mu}^{\sigma} f := \sigma^2 \Delta _\MM f + \left( \nabla _\MM (\sigma^2) - \sigma^2 \nabla_\MM V \right) \nabla _\MM f,
$$
where $\Delta _\MM$ is the Laplace-Beltrami operator on $\MM$ and $\sigma^2$ is the weight function appearing in the weighted Poincar\'e inequality, which now rewrites by integration by parts as the inequality
$$
\lambda \Var_\mu (f) \leq - \int_\MM f \LL _\mu ^\sigma f d\mu .
$$
When $\sigma \equiv 1$, the operator reduces to the classical dynamics
$$
\LL_{\mu} f := \Delta _\MM f - \nabla_\MM V \nabla_\MM f.
$$
Therefore the best constant in the weighted Poincar\'e inequality \eqref{modipoin} is nothing but the spectral gap in $L^2(\mu)$ of the operator $- \LL_\mu ^\sigma$ (or more rigorously the self-adjoint extension of this operator, still denoted $- \LL_\mu ^\sigma$ in the remainder of the paper). Recall that the spectral gap of the operator $- \LL_\mu ^\sigma$ is, when it exists, the largest real $\lambda$ such that its spectrum lies in $\{ 0 \} \cup [\lambda, +\infty )$. We denote it $\lambda_1 (-\LL _\mu ^\sigma)$. It is characterized by the variational identity
\begin{equation}
\label{eq:varia}
\lambda_1 (-\LL _\mu ^\sigma) = \inf \left \{ \frac{- \int_\MM f \LL _\mu ^\sigma f d\mu}{\Var_\mu (f)} : f \in H^1 _\sigma (\MM,\mu) , f \neq \mbox{const} \right \} ,
\end{equation}
where $H^1 _\sigma (\MM , \mu) $ is the weighted Sobolev space
$$
H^1 _\sigma (\MM,\mu) := \left \{ f\in L^2 (\mu) : \int_\MM \sigma^2 \vert \nabla _\MM f \vert ^2 d\mu < \infty \right \} ,
$$
which is nothing but the domain of the associated Dirichlet form. According to this point of view, the difficulty to establish a weighted Poincar\'e inequality with weight $\sigma^2$ and to estimate the optimal constant is transferred to the study of the spectral gap of the operator $-\LL _\mu ^\sigma$ which can be more tractable, in particular in the one-dimensional setting, in view of the recent results emphasized in \cite{bj}. \vspace{0.1cm}

From now on we turn to the one-dimensional framework of the real line, for which the diffusion operator is
$$
\LL_{\mu}^{\sigma} f := \sigma^2 f'' + \left( (\sigma^2)' - \sigma^2 V' \right) f' .
$$
To deal with an essentially self-adjoint operator, we require the following set of assumptions on $\sigma$, which will be satisfied by our examples all along the paper:

$\circ$ Smoothness: the function $\sigma \in \CC ^\infty (\rr)$.

$\circ$ Ellipticity: the function $\sigma$ is positive on $\rr$.

$\circ$ Completeness: the metric space $(\rr,d_\sigma)$ is complete, where $d_\sigma$ denotes the distance $d_\sigma (x,y) = \left \vert \int_x ^y du /\sigma(u)\right \vert $. \vspace{0.1cm}

Beyond the famous Muckenhoupt estimates \cite{muckenhoupt} which can be limited in some cases of interest, a practical criterion ensuring a spectral gap for the dynamics has been put forward in \cite{bj} by the two first authors. Using the intertwining approach, one of the main results can be reformulated as follows:
\begin{equation}
\label{e:chen}
\lambda_1 (-\LL_\mu ^\sigma) \geq \sup_{f\in \II \CC ^\infty (\rr)} \inf  V_f ^\sigma ,
\end{equation}
where $\II \CC^\infty (\rr)$ is the set of functions in $\CC^\infty (\rr)$ with positive derivative and
$$
 V_f ^\sigma := \frac{(- \LL _\mu ^\sigma f)'} {f'} .
$$
Moreover the equality holds trivially in \eqref{e:chen} if $\lambda_1 (-\LL_\mu ^\sigma)$ is an eigenvalue. Such a formula, which is reminiscent of the work of Chen and Wang \cite{chen_wang} (see also \cite{Chen99} and \cite{djellout} through various alternative forms), provides a convenient criterion ensuring a spectral gap for the dynamics, at least when one has the intuition of the shape of the associated eigenfunctions. For instance if one choose the family of functions $f_\varepsilon ' = e^{\varepsilon V} / \sigma^2$ indexed by some convenient real number $\varepsilon$, then after some computations one obtains $\LL_\mu ^\sigma f = -(1-\varepsilon) V' e^{\varepsilon V}$ and
\begin{equation}
\label{eq:v_f}
V_{f_\varepsilon} ^\sigma = (1-\varepsilon) \sigma^2 \left( V'' + \varepsilon {V'}^2 \right) ,
\end{equation}
leading to a potentially interesting lower bound on the spectral gap. \vspace{0.1cm}

Recently, an estimate of the type \eqref{e:chen} has been successively used by the authors in \cite{bjm} to derive classical and weighted Poincar\'e inequalities for spherically symmetric probability measures (including the log-concave and heavy tailed cases) in dimension larger than 2, after a reduction of the problem to the one-dimensional setting. However the one-dimensional case, for which an additional analysis is required, has not been treated yet and this is what we intend to present in this short note. In particular we focus our attention on the one-dimensional version of the examples studied in \cite{bjm}. \vspace{0.1cm}

This work is presented as follows. In Section \ref{sect:power}, we investigate a first family of examples which already appears in \cite{bj} (as an illustration of the intertwining techniques for deriving functional inequalities) and for which we establish classical and weighted Poincar\'e inequalities. The measure $\mu$ is the so-called exponential power distribution of parameter $\alpha \in [1,\infty)$, that is to say, the potential $V$ is convex and is given by $V(x) = |x|^{\alpha}/\alpha$. Using first \eqref{e:chen} in a convenient way, we are able to derive new estimates for the spectral gap $\lambda_1 (-\LL _\mu)$ which are different as $\alpha \in [1,2]$ or $\alpha > 2$, showing the criticality of the Gaussian case $\alpha =2$. Next, to illustrate the robustness of formula \eqref{e:chen}, we show in this Gaussian setting how to obtain a family of weighted Poincar\'e inequalities that refine the classical one, the weight function being $\sigma ^2 (x) = 1/(1+ bx^2)$ where $b$ is some positive parameter. In particular we are able to capture the optimal constant $\lambda_1 (-\LL ^\sigma _\mu)$, which is different as $b \in (0,1/2)$ or $b \geq 1/2$ but continuous at point $b = 1/2$. We pursue this approach in Section \ref{sect:cauchy} where we investigate the case of the generalized Cauchy distribution with potential $V(x) = \beta \log (1+ x^2)$, $\beta >1/2$. Because of the heavy-tailed phenomenon, such a probability measure does not satisfy the classical Poincar\'e inequality but a (weaker) weighted Poincar\'e inequality involving the weight $\sigma ^2 (x) = 1+ x^2$. Once again the optimal constant $\lambda_1 (-\LL ^\sigma _\mu)$ is obtained and it reveals to be different as $\beta \in (1/2 ,3/2]$ or $\beta >3/2$, with also a continuous transition as point $\beta = 3/2$.

\section{Classical and weighted Poincar\'e inequalities for the exponential power distributions}
\label{sect:power}
The first class of examples we consider in this part is of the following form: the potential is $V(x) = |x|^{\alpha}/\alpha$ with $\alpha \geq 1$, so that the underlying probability measure, the exponential power distribution of parameter $\alpha$, has Lebesgue density on $\rr$ given by
$$
\mu (dx)= \frac{1}{2 \alpha^{1/\alpha -1} \Gamma(1/\alpha)} e^{-|x|^{\alpha}/\alpha} dx ,
$$
where $\Gamma$ is the famous Gamma function, $\Gamma (z) = \int_0 ^\infty t^{z-1} e^{-t} dt$, $z>0$. Let us start by the case of classical Poincar\'e inequalities.

\subsection{Classical Poincar\'e inequalities}
\label{sect:classic}
By the convexity of the potential $V$, the measure $\mu$ satisfies a classical Poincar\'e or, in other words, the associated classical Markovian operator given for any $f\in \CC _0 ^\infty (\rr)$ by
$$
\LL_{_\mu} f(x) = f''(x)-|x|^{\alpha-2} x f'(x),
$$
has a spectral gap. Among the potential cases of interest, some are well-studied in the literature:
\begin{itemize}
\item[$\circ$] when $\alpha=1,$ the measure $\mu$ is the double-exponential distribution and we have $\lambda_1(-\LL_\mu)=1/4;$
\item[$\circ$] when $\alpha=2,$ we consider the standard Gaussian distribution for which we have $\lambda_1(-\LL_\mu)=1;$
\item[$\circ$] when $\alpha = \infty$, this is the uniform measure on $[-1, 1]$ and we have $\lambda_1(-\LL_{\mu})=\pi^2/4.$
\end{itemize}
Actually, it seems difficult to obtain the exact value of the spectral gap beyond these three cases. Below, we propose new estimates which are rather different as $\alpha \in [1,2]$ or $\alpha > 2$. Such an observation is expected because the Gaussian case $\alpha =2$ is critical in terms of convexity of $V$: the infimum of the second derivative of $V$ is attained at infinity when $\alpha \in [1,2)$ and at the origin when $\alpha > 2$. In particular, our estimates are sharp as $\alpha$ is close to 1 or 2. However there is still room for a mild improvement in the extremal regime $\alpha \to \infty$ since we do not capture the $\pi^2/4$ ($\approx 2.467$) by a limiting argument (we only obtain the interval $[2,3]$ at the limit $\alpha \to \infty$). The result is stated as follows.
\bthm
\label{main} If $\mu$ stands for the exponential power distribution of parameter $\alpha \in [1,\infty)$, then the spectral gap $\lambda_1(-\LL_{\mu})$ of the classical Markovian operator $-\LL_{\mu}$ satisfies the following estimates:
\begin{itemize}
\item[$(i)$] if $1\le \alpha \leq 2$, we get
$$
\frac{\alpha^2}{4}\le\lambda_1(-\LL_\mu) \le 2^{1 - 2/\alpha };
$$
\item[$(ii)$] if $\alpha\ge 2,$ then we have
$$
\frac{2(1+\alpha)^{1-2/\alpha}}{\alpha} \leq \lambda_1(-\LL_\mu)\leq 3^{1-2/\alpha} .
$$
\end{itemize}
\nthm

\bprf
The proof is divided in several parts. Let us first concentrate our attention on the two upper bounds. By the variational formula \eqref{eq:varia}, a brief computation with the function $f(x)=x$ entails
$$
\lambda_1(-\LL_{\mu})\leq \alpha^{-2/\alpha}\frac{\Gamma(1/\alpha)}{\Gamma(3/\alpha)}.
$$
Now our objective is to understand the behaviour of the right-hand-side in terms of the parameter $\alpha$ and we thus have to distinguish the two cases: \vspace{0.1cm}

$\circ$ On the one hand if $\alpha \geq 2$ then our idea is to use the log-convexity of the Gamma function. More precisely, if $a,b$ are two parameters such that $a>0$ and $b\in [0,1]$, then the log-convexity of $\Gamma$ entails the inequality
$$
a \Gamma (a) = \Gamma (a+1) \leq \Gamma (a+b) ^{b} \, \Gamma (a+b+1) ^{1-b} = (a+b)^{1-b} \, \Gamma (a+b) ,
$$
which rewrites in a condensed form as
$$
\frac{\Gamma (a) \, a^b}{\Gamma (a+b)} \leq \left( \frac{a+b}{a} \right) ^{1-b} .
$$
Using then this estimate with $a = 1/ \alpha $ and $b = 2/ \alpha $ yields the result. \vspace{0.1cm}

$\circ$ On the other hand if $\alpha \in [1,2]$, we proceed differently and use the famous Kershaw inequality \cite{kershaw} controlling the ratio of two Gamma functions: for any $(a,b) \in (0,\infty) \times [1,\infty)$ such that $a \leq b \leq a+1$,
$$
\frac{\Gamma (a)}{\Gamma (b)} \leq \left(\frac{a+b-1}{2}\right) ^{a-b} .
$$
Then the desired upper estimate holds with the two parameters $a= 1+1/\alpha $ and $b = 3/\alpha $, which verify the required conditions above. \vspace{0.1cm} \\
The second step of the proof is to establish the lower bound in the case $\alpha \in [1,2]$. To that aim, we would like to use the inequality \eqref{e:chen}. Note that the measure $\mu$ interpolates between the double-exponential and Gaussian distributions. In particular, the spectral gap is attained in the Gaussian case for linear functions whereas for the double-exponential distribution we obtain its exact value by choosing in the variational formula \eqref{eq:varia} the family of functions $f_a (x) = e^{a \vert x\vert }$ and letting $a \uparrow 1/2$. Hence our idea is to try in \eqref{e:chen} some functions of the type $f_\varepsilon '(x) = e^ { \varepsilon\vert x\vert ^{\alpha} /\alpha}$, where $\varepsilon $ is some constant depending on $\alpha$ and which equals $1/2$ when $\alpha =1$ and vanishes when $\alpha = 2$. Note that such a choice of functions $f_\varepsilon '$ is of the form appearing in \eqref{eq:v_f} (with $\sigma \equiv 1$) and therefore we have
$$
\LL_\mu f_\varepsilon  (x) = -(1-\varepsilon) \vert x\vert ^{\alpha-2} x e^{\varepsilon \vert x\vert ^\alpha /\alpha} ,
$$
and
$$
V_{f_\varepsilon } (x) = \frac{(-\LL_\mu f_\varepsilon )'(x)}{f_\varepsilon  '(x)} = (1-\varepsilon) \left( (\alpha-1)\vert x\vert ^{\alpha-2} + \varepsilon \vert x\vert ^{2(\alpha-1)} \right),
$$
so that the infimum is positive as soon as $\varepsilon \in (0,1)$ and thus equals
$$
\inf V_{f_\varepsilon } =  (1-\varepsilon) \left( (\alpha-1) \left( \frac{2-\alpha}{2\varepsilon}\right) ^{1-2/\alpha} + \varepsilon \left( \frac{2-\alpha}{2\varepsilon}\right) ^{2-2/\alpha} \right).
$$
Finally choosing $\varepsilon = (2-\alpha) /2$, we obtain through \eqref{e:chen} the announced lower bound $\lambda_1(-\LL_{\mu})\geq \alpha ^2 /4$. \vspace{0.1cm} \\
Let us achieve this proof by showing the lower bound when $\alpha \geq 2$. As in the previous case, we want to use the estimate \eqref{e:chen}. Since the choice of functions $f_\varepsilon '$ involved in the identity \eqref{eq:v_f} does not work (we have no idea of the underlying eigenfunctions), a reasonable strategy is to find a convenient centered function $\rho \in \II \CC ^\infty (\rr)$ such that the solution $f$ to the Poisson equation $-\LL _\mu f = \rho$ also belongs to the set $\II \CC ^\infty (\rr)$. In other words, the key point is to understand the action of the inverse operator $(-\LL_\mu)^{-1}$ on a class of smooth increasing functions. Certainly, such an operator is not well-defined on $\CC ^\infty (\rr)$ but (sometimes) on a convenient subspace of smooth centered Lipschitz functions, the metric under consideration being induced by $\rho$, cf. \cite{djellout} for more details.
Hence, let us start by the simplest centered function $\rho \in \II \CC ^\infty (\rr)$ we may choose, namely the identity function $\rho(x) = x$. Solving then the Poisson equation (a differential equation of the first order in the derivative $f'$) yields
$$
f'(x) = e^{|x|^{\alpha}/\alpha} \int_{x}^{\infty} y e^{-|y|^{\alpha}/\alpha}dy.
$$
Since it belongs to the set $\II \CC ^\infty (\rr)$, it is a good candidate for plugging into the formula \eqref{e:chen}. We have
$$
\inf _{x\in \rr } V_f (x) = \frac{1}{\sup_{x\in \rr} f'(x)} = \frac{1}{\sup_{x\geq 0}\alpha^{2/\alpha-1}e^x\int_x^{\infty}t^{2/\alpha-1} e^{-t} dt },
$$
where we used a symmetry argument and a change of variable in the integral to obtain the last equality. Differentiating the denominator with respect to $x$ shows that the supremum over $[0,\infty)$ is attained at the origin since $\alpha$ lies in the region $[2,\infty)$ (it would be infinite otherwise) and we thus get the lower bound
$$
\lambda_1(-\LL_{\mu}) \geq \frac{\alpha^{1 -2/\alpha}}{\Gamma(2/\alpha)} .
$$
Finally, we apply once again Kershaw's estimate with $a = 1+ 1/\alpha $ and $b = 2$ to get the required estimate $\lambda_1(-\LL_\mu) \geq 2\alpha ^{-1} (1+\alpha)^{1-2/\alpha} $. The proof is now complete.
\nprf

\subsection{Weighted Poincar\'e inequalities for the Gaussian distribution}

As mentioned above, we illustrate now how the formula \eqref{e:chen} might be used to get a weighted Poincar\'e inequality instead of a classical one, as stated in the last section. For more simplicity, let us focus our attention on the Gaussian setting, that is, the parameter $\alpha$ equals 2, although the forthcoming analysis might be adapted to any $\alpha \geq 1$ (the same approach might be applied and one can prove that the weighted Poincar\'e inequality with (optimal) weight $\sigma^2 (x)= (1+x^2)^{1-\alpha}$ holds). The weight function we consider is of the type $\sigma_{a,b} ^2(x)=1/(a+ b x^2)$, with $a,b>0$, and the Markovian dynamics is now given by
\begin{eqnarray*}
\LL_{\mu}^{\sigma_{a,b}} f(x) & = & \sigma_{a,b} ^2 (x) f''(x) + \left( (\sigma_{a,b} ^2)'(x) - \sigma_{a,b} ^2 (x) V'(x) \right) f'(x) \\
& = & \frac{f''(x)}{a+bx^2} - \frac{(2b+a)x + bx^3}{(a+bx^2)^2}  f'(x) .
\end{eqnarray*}
We will see later why such a weight is optimal. Since we have clearly $\lambda_1 (-\LL_\mu ^{\sigma_{a,b}}) = a \lambda_1 (-\LL_\mu ^{\sigma_{1,b/a}})$, we can assume without loss of generality that the weight is given in the sequel by $\sigma^2 (x) := \sigma_{1,b} ^2 (x)=1/(1+ b x^2)$. Below, we obtain a family (indexed by $b$) of weighted Poincar\'e inequalities that strengthen the classical Poincar\'e inequality for the Gaussian measure. In particular, an interesting feature resides in the continuous transition between two different regimes as $b$ lies in the regions $(0,1/2 )$ or $[1/2, \infty)$, the reason being a lack of integrability in the second case. Such a phenomenon is also shared by an heavy-tailed probability measure such as the generalized Cauchy distribution, as we will observe in Section \ref{sect:cauchy}. For the moment, our estimates in the Gaussian setting stand as follows.
\bthm
\label{thmforgauss} If $\mu$ denotes the standard Gaussian measure and $\sigma$ is the weight $\sigma^2 (x) =1/(1+ b x^2)$, then the spectral gap of the operator $-\LL^{\sigma}_{\mu}$ is given by
\[
\lambda_1 (-\LL_\mu ^\sigma) = \left \{
\begin{array}{lll} 1-b & \mbox{if} & 0 <b < 1/2; \vspace{0.1cm} \\
1/4b & \mbox{if} & b \geq 1/2. \\
\end{array}
\right.
\]
In other words, the spectral gap given above is the best constant $\lambda >0$ such that we have the following weighted Poincar\'e inequality: for any $f\in \CC ^\infty _0 (\rr)$,
$$
\lambda \Var _{\mu}(f) \leq \int_{\rr} \frac{f'(x)^2}{1+b x^2}\mu(dx) .
$$
\nthm
\bprf
If $b\in (0,1/2)$ then with the choice of the centered function $f(x)=xe^{bx^2 /2}$, we have $\LL_\mu ^{\sigma} f = - (1-b) f$ and since $f$ is increasing and in $L^2(\mu)$ (because $0 < b <1/2$), the value $1-b $ is nothing but the expected spectral gap. \vspace{0.1cm} \\
Now let us focus our attention on the region $b \geq 1/2$ for which the latter argument is no longer available. In contrast to Section \ref{sect:classic} where there is no weight, the contribution of the weight $\sigma^2$ in the generator has to be taken into account and a natural candidate to invoke the formula \eqref{e:chen} is to pick a family of functions according to the formula \eqref{eq:v_f}, i.e. of the type $f_\varepsilon '(x) = (1+bx^2) e^{\varepsilon x^2 /2}$ and then to check the best $\varepsilon$ allowing us to use \eqref{e:chen}. We obtain
$$
\LL_\mu ^\sigma f_\varepsilon  (x) = -(1-\varepsilon) x e^{\varepsilon x^2 /2} \quad \mbox{ and } \quad V_{f_\varepsilon } ^\sigma (x) = \frac{(-\LL_\mu ^\sigma f_\varepsilon )'(x)}{f_\varepsilon  '(x)} = (1-\varepsilon) \frac{1+\varepsilon x^2 }{1+bx^2} ,
$$
and provided $\varepsilon \in (0,1)$ we have
\[
\inf V_{f_\varepsilon } ^\sigma =
\left\{
\begin{array}{lll}
\varepsilon (1-\varepsilon)/b & \mbox{if} & 0 < \varepsilon \leq b ; \vspace{0.1cm} \\
(1-\varepsilon) & \mbox{if} &  b \leq \varepsilon <1 .
\end{array}
\right.
\]
Finally optimizing in $\varepsilon \in (0,1)$ yields the lower bound $\lambda_1(-\LL^{\sigma}_{\mu}) \geq 1/4b$. \vspace{0.1cm} \\
To establish the upper bound, we choose the function $f(x)=x e^{(1-\varepsilon)x^2/4}$ with $\varepsilon \in (0,1)$ to guarantee the existence of the forthcoming quantities. We have
$$\Var _\mu (f) = \frac{1}{\sqrt{2\pi}}\int_{\rr} y^2 e^{-\varepsilon y^2/2}dy = \frac{1}{\varepsilon^{3/2} },
$$
and
\begin{eqnarray*}
- \int_\rr f \LL_\mu ^\sigma f d\mu & = &  \int_\rr \sigma^2 f'^2 d\mu \\
& = & \int_{\rr} \left( 1+\frac{(1-\vep ) y^2}{2} \right)^2 \frac{e^{-\vep y^2/2}}{1+by^2} \frac{dy}{\sqrt{2\pi}} \\
& = & \int_{\rr} \left( \frac{(1-\varepsilon)^2 y^4}{4} +(1-\vep ) y^2 + 1\right) \frac{e^{-\vep y^2/2}}{1+by^2} \frac{dy}{\sqrt{2\pi}} \\
& \leq & \frac{(1-\varepsilon)^2}{4b} \int_{\rr} y^2 e^{-\vep y^2/2} \frac{dy}{\sqrt{2\pi}} + \frac{1-\varepsilon +b }{b} \int_{\rr} e^{-\vep y^2/2} \frac{dy}{\sqrt{2\pi}} \\
& = & \frac{(1-\varepsilon)^2}{4b\varepsilon^{3/2}} + \frac{1-\varepsilon + b}{b\sqrt{\varepsilon}} ,
\end{eqnarray*}
so that by the variational identity \eqref{eq:varia}, we get
\begin{eqnarray*}
\lambda_1(-\LL_\mu ^{\sigma}) & \leq & \frac{- \int_\rr f \LL_\mu ^\sigma f d\mu}{\Var _\mu (f) } \\
& \leq & \frac{(1-\varepsilon)^2}{4b} + \frac{(1-\varepsilon + b)\varepsilon }{b} .
\end{eqnarray*}
Finally letting $\vep$ tend to $0$ guarantees the desired upper bound. The proof of the equality for $b\geq 1/2$ is now achieved.
\nprf
Note that the lower bound $1/4b$ in the case $b\geq 1/2$ can be derived directly from the case $b=1/2$. Indeed using the variational formula \eqref{eq:varia} of the spectral gap shows that we have for any $b \geq 1/2$ the comparison $\lambda_1 (-\LL_\mu ^{\sigma_{1,b}}) \geq \lambda_1 (-\LL_\mu ^{\sigma_{1,1/2}}) /2b$. \vspace{0.1cm}

As announced at the beginning of this part, the polynomial of degree 2 at the denominator of the weight $\sigma$ is optimal among all positive polynomials, and also among all powers of $x$. More precisely, if $\sigma$ is of the form $\sigma ^2 (x) = (1+x^2)^{-a}$, where $a$ is some parameter very close to 1 (but greater than 1), then we have $\lambda_1 (-\LL_\mu ^\sigma) = 0$. Indeed, we have with the functions $f_\varepsilon (x) = e^{\varepsilon x^2 /2}$, where $\varepsilon \in (0,1/2)$,
$$
\Var_\mu (f_\varepsilon ) = \frac{1}{\sqrt{1-2\varepsilon}} - \frac{1}{1-\varepsilon} ,
$$
and also
$$
- \int_\rr f_\varepsilon \LL_\mu ^\sigma f_\varepsilon  d\mu = \frac{2^{2-a} \varepsilon ^2 \Gamma (3/2 - a )}{\sqrt{2\pi} (1-2\varepsilon )^{3/2 - a}} ,
$$
so that plugging these two quantities into the variational identity \eqref{eq:varia} and taking then the limit as $\varepsilon \uparrow 1/2$ yields the conclusion. \vspace{0.1cm}

Before turning to generalized Cauchy distributions in the next section, let us say some words about another functional inequality related to weighted Poincar\'e inequalities, known as the Brascamp-Lieb inequality \cite{brascamp-lieb}. For a given measure $\mu$ with Lebesgue density proportional to $e^{-V}$, we assume that the potential $V$ is strictly convex on $\rr$ and moreover $V' \in L^2(\mu)$. Such a log-concave framework includes for instance the case of the exponential power distribution of parameter $\alpha$, with $\alpha$ in the region $[1,2]$. Then the Brascamp-Lieb inequality reads as follows: for any $f\in \CC ^\infty _0 (\rr)$, we have
$$
\Var_\mu (f) \leq \int_\rr \frac{f'^2}{V''} d\mu ,
$$
and the constant 1 in front of the variance is optimal. Let us recover briefly this inequality by a simple argument. Note that the Brascamp-Lieb inequality can be rewritten as a weighted Poincar\'e inequality with the weight function $\sigma^2=1/V''$. Therefore the best constant in this inequality is nothing but the spectral gap of the Markovian dynamics
\begin{eqnarray*}
\LL_{\mu}^{\sigma} f & = & \sigma^2 f'' + \left( (\sigma^2)' - \sigma^2 V' \right) f' \\
& = & \frac{f''}{V''} - \left( \frac{V'''}{{V''}^2} + \frac{V'}{V''} \right) f' .
\end{eqnarray*}
Once we have in mind this analogy, we observe that the equality $\LL_{\mu}^{\sigma} V' = -V'$ holds, i.e., the function $V'$ is always an eigenfunction associated to the eigenvalue 1 (recall that $V' \in L^2(\mu)$) and since $V'$ is increasing, we deduce that $\lambda_1 (-\LL_\mu ^\sigma) = 1$. However this Brascamp-Lieb inequality is somewhat limited in the sense that the weight function is fixed in terms of the potential $V$, excluding \textit{de facto} a large variety of interesting weighted inequalities (for instance in the Gaussian case we have $V'' \equiv 1$, leading only to the classical Poincar\'e inequality). Recently, the two first authors obtained a generalization of the Brascamp-Lieb inequality, involving the function $V_f$ appearing in the formula \eqref{e:chen}. It is stated as follows. Consider the classical operator $\LL_\mu f = f'' - V' f'$, where $V$ is some smooth potential, non necessarily convex. If there exists $f\in \II \CC ^\infty (\rr)$ such that the function $V_f$ defined in \eqref{e:chen} is positive, then we have for any $g\in \CC ^\infty _0 (\rr)$,
$$
\Var_\mu (g) \leq \int_\rr \frac{g'^2}{V_f} d\mu .
$$
When $V$ is strictly convex, we recover the Brascamp-Lieb inequality above by considering the identity function $f(x) = x$. However choosing conveniently the function $f$ allows a large variety of weighted Poincar\'e inequalities. For instance consider in \eqref{e:chen} the function $f' = e^{-W}$ where $W$ has to be carefully chosen. After some brief computations, one has
$$
V_f = W'' - {W'}^2 - V' W' +V'' ,
$$
and rewriting this expression with the new function $U := W' + V' /2$, we get a more tractable function
$$
V_f = U' - U^2 + \frac{{V'}^2}{4} + \frac{V''}{2} .
$$
Now the main point is to choose the function $U$ in such way that $V_f$ is positive. To that aim, we can consider for instance a polynomial and the easiest form one can think of is $U(x) = \gamma x$, where $\gamma$ is some real parameter. The function of interest is then $f'(x) = e^{(V(x) - \gamma x^2) /2}$ and we get
$$
V_f (x) = \gamma - \gamma ^2 x^2 + \frac{{V'(x) }^2}{4} + \frac{V''(x)}{2} ,
$$
which can be positive for some given potentials $V$ (provided $\gamma$ is suitably chosen). In particular such an argument might be used when the function $V'^2 + 2V''$ is not positive everywhere. In terms of weighted Poincar\'e inequalities, it leads to the inequality
$$
\frac{1}{4} \Var _\mu (g) \leq \int_\rr \frac{{g'(x) }^2 }{4 \gamma - 4 \gamma ^2 x^2 + {V'(x) }^2 + 2 V''(x)} \mu (dx).
$$
Coming back to the Gaussian case $V(x) = x^2 /2$, we have
$$
V_f (x) = \left( \frac{1}{4}- \gamma ^2 \right) x^2 + \gamma + \frac{1}{2},
$$
which is positive for any $\gamma \in [0,1/2]$. Therefore we obtain the family (indexed by $\gamma$) of inequalities
$$
\left( \gamma + \frac{1}{2} \right) \Var _\mu (g) \leq \int_\rr \frac{{g'(x) }^2 }{1+ (1/2 - \gamma ) x^2 } \mu (dx).
$$
These inequalities, which reduce as expected to the classical Poincar\'e inequality when $\gamma = 1/2$, enter the general framework of Theorem \ref{thmforgauss} with $b = 1/2 - \gamma \in [0,1/2]$, for which we recover without any additional effort the optimal constant $1 - b = \gamma +1/2$.

\section{Weighted Poincar\'e inequalities for generalized Cauchy distributions}
\label{sect:cauchy}
This final part is devoted to study the case of an heavy-tailed distribution, namely the generalized Cauchy distribution of parameter $\beta >1/2$, that is,
$$
\mu (dx)=\frac{dx}{Z(1+x^2)^{\beta}} , \quad x\in \rr ,
$$
where $Z$ is the devoted normalization constant. In the language of the previous sections, the potential $V$ rewrites as $V(x) = \beta \log (1+x^2)$ and therefore is far from being convex at infinity. In particular, the classical Markovian operator $\LL_\mu f = f'' - V'f'$ doe not have a spectral gap or, in other words, $\lambda_1( - \LL_\mu)=0$. In terms of functional inequality, the measure does not satisfy the classical Poincar\'e inequality, that is, there is no constant $\lambda >0$ such that we have
$$
\lambda \Var _\mu (f) \leq \int_\rr  f'^2 d\mu , \quad f\in \CC ^\infty _0 (\rr).
$$
However, one may hope to get a weighted Poincar\'e inequality that is weaker than the classical one, and this is indeed the case by choosing the weight function $\sigma ^2 (x) = 1+x^2$, cf. \cite{BLe,bj}. To that aim, we consider the Markovian dynamics
\begin{eqnarray*}
\LL_{\mu}^{\sigma} f(x) & = & \sigma^2 (x) f''(x) + \left( (\sigma^2)'(x) - \sigma^2 (x) V'(x) \right) f'(x) \\
& = & (1+x^2) f''(x) + 2(1-\beta)x f'(x) .
\end{eqnarray*}
Below, we give the exact value of the spectral gap $\lambda_1 (-\LL_\mu ^\sigma)$ of the operator $-\LL_\mu ^\sigma$. In particular, our results exhibit as in the Gaussian case above a continuous transition at point $\beta = 3/2$, which is somewhat expected due to a lack of integrability for $\beta$ in the region $(1/2 , 3/2)$.
\bthm
\label{thmforcauchy}
Let $\mu$ be the generalized Cauchy distribution of parameter $\beta >1/2$ and let $\sigma$ be the weight function $\sigma ^2 (x) = 1+x^2$. Then the spectral gap $\lambda_1 (-\LL_\mu ^\sigma)$ of the operator $-\LL_\mu ^\sigma$ satisfies:
\[
\lambda_1 (-\LL_\mu ^\sigma) = \left \{
\begin{array}{lll} 2(\beta-1) & \mbox{if} & \beta > 3/2; \vspace{0.1cm} \\
\left(\beta-1/2 \right)^2  & \mbox{if} & 1/2 < \beta \leq 3/2. \\
\end{array}
\right.
\]
In other words, the spectral gap given above is the best constant $\lambda >0$ such that we have the following weighted Poincar\'e inequality: for any $f\in \CC ^\infty _0 (\rr)$,
$$
\lambda \Var_\mu (f) \leq \int_{\rr} (1+x^2) (f'(x))^2 \mu (dx).
$$
\nthm

\bprf
As in the Gaussian setting, we start by the case for which the eigenfunction is easily found, that is for $\beta>3/2$. Indeed taking the centered function $f(x)=x$, we have $\LL_\mu ^{\sigma} f = -2 (\beta-1) f$ and since $f$ is increasing and lies in $L^2(\mu)$, the value $2 (\beta-1)$ is nothing but the spectral gap. \vspace{0.1cm} \\
When $\beta \in (1/2 , 3/2]$, we proceed for the lower bound as in the Gaussian case studied previously and choose according to \eqref{eq:v_f} the functions $f_\varepsilon '(x) = (1+x^2) ^{\varepsilon}$ for some real parameter $\varepsilon$ (the present $\varepsilon$ corresponds to $ \varepsilon \beta -1$ in the notation of \eqref{eq:v_f}). Then we have $\LL_\mu ^\sigma f(x) = -2 (\beta -\varepsilon -1)x(1+x^2)^\varepsilon$ and also
$$
V_{f_\varepsilon} (x) = \frac{(- \LL_\mu ^\sigma f_\varepsilon)'(x)}{f_\varepsilon ' (x)} = \frac{2 (\beta -\varepsilon - 1)(1+(2\varepsilon+1)x^2)}{1+x^2} .
$$
Now if $\varepsilon$ is chosen in the interval $(-1/2 , \beta-1)$ then $\inf V_{f_\varepsilon} = 2 (\beta -\varepsilon - 1)(2\varepsilon+1) >0$ and optimizing then in $\varepsilon$ (the optimal one is $\varepsilon = (2\beta -3)/4 \in (-1/2 , \beta-1)$) yields through the formula \eqref{e:chen} the lower estimate $\lambda_1(-\LL_{\mu}^{\sigma})\geq (\beta-1/2)^2$. \vspace{0.1cm} \\
To obtain the upper bound, take the function $f(x)=x(1+x^2)^{\varepsilon}$ with $\varepsilon < (2\beta -3)/4.$ After some careful calculations and using the notation $a := \beta-2\varepsilon >3/2$, we have
$$
\Var_\mu (f)= \frac{\Gamma(a- 3/2)\Gamma(1/2)}{Z\Gamma(a-1)}\bigg(1-\frac{a-3/2}{a-1}\bigg),
$$
and
$$
- \int_\rr f\LL_\mu ^\sigma f d\mu = \frac{\Gamma(a- 3/2)\Gamma(1/2)}{Z\Gamma(a-1)} \bigg((1+2\varepsilon)^2-4\varepsilon(1+2\varepsilon)\frac{a-3/2}{a-1}+4\varepsilon^2\frac{(a-3/2)(a-1/2)}{a(a-1)}\bigg).$$
Using then the variational formula \eqref{eq:varia} entails
\begin{eqnarray*}
\label{uppcauchy}
\lambda_1(-\LL_{\mu}^{\sigma}) &\leq & 2\bigg((1+2\varepsilon)^2 (a-1)-4\varepsilon(1+2\varepsilon)(a-3/2)+4\varepsilon^2\frac{(a-3/2)(a-1/2)}{a}\bigg) \\
& = & 2\left( a-1+2\varepsilon+\frac{3\varepsilon^2}{a} \right) \\
& = & 2(\beta-1)+\frac{6\varepsilon^2}{\beta-2\varepsilon}.
\end{eqnarray*}
Finally letting $\varepsilon \uparrow (2\beta-3)/4$, we get
$$\lambda_1(-\LL_{\mu}^{\sigma}) \leq 2(\beta-1)+ \frac{6(2\beta-3)^2}{16(\beta-2(2\beta-3)/4)} = \left( \beta-\frac12 \right) ^2,
$$
which is the required upper bound. The proof is now complete.
\nprf

Let us conclude this work by two remarks. First, we mention that similarly to the Gaussian case, the weight in the weighted Poincar\'e inequality is optimal among all powers of $x$, that is, if the weight function is chosen to be $\sigma^2 (x)=(1+x^2)^b$ with some $b\in (0,1)$, then we have $\lambda_1(-\LL_{\mu}^{\sigma }) =0.$ Indeed, choose the function  $f(r)=(1+r^2)^{\vep},$ with $0<\vep<\beta /2 - 1/4.$ Then after some tedious computations, we have with the notation $a := \beta -2\varepsilon >1/2$,
$$
\mu(f) = \frac{\Gamma(\beta-\vep- 1/2)\Gamma(\beta)}{\Gamma(\beta- 1/2)\Gamma(\beta-\vep)} , \quad \mu(f^2)=\frac{\Gamma(a-1/2)\Gamma(\beta)}{\Gamma(\beta- 1/2)\Gamma(a)} ,
$$
and
$$
- \int_\rr f \LL _\mu ^{\sigma } f d\mu = 2 \vep^2\frac{\Gamma(a+1/2-b)\Gamma(\beta)}{\Gamma(\beta- 1/2)\Gamma(a+2-b)}.
$$
Therefore we obtain by the variational formula \eqref{eq:varia}
$$
\lambda_1(-\LL_{\mu}^{\sigma })^{-1}\ge \bigg(\frac{\Gamma(a- 1/2)}{\Gamma(a)}-\frac{\Gamma^2(\beta-\vep- 1/2)\Gamma(\beta)}{\Gamma(\beta- 1/2)\Gamma^2(\beta-\vep)}\bigg)\frac{\Gamma(a+2-b)}{2\vep^2\Gamma(a+1/2-b)}.$$
When $\vep$ tends to $\beta/2-1/4,$ the limit of the right-hand-side of the above inequality is infinite or, in other words, $\lambda_1(-\LL_{\mu}^{\sigma }) = 0$. \vspace{0.1cm}

Our second brief observation follows an argument due to Bobkov and Ledoux \cite{BLe}. Given $\beta >1/2$, we slightly modify the measure $\mu$ and replace it  by
$$
\mu _\beta(dx)= \left( 1+\frac{x^2}{2\beta-1}\right) ^{-\beta} \frac{dx}{Z}.
$$
Then it is clear that $\mu_\beta$ converges weakly as $\beta$ tends to infinity to the standard Gaussian measure. Since $\beta$ has the vocation to be large, we assume $\beta >3/2$ and according to the previous results, we have
$$
\lambda_1 (-\LL_{\mu_\beta} ^{\sigma }) = \frac{2(\beta-1)}{2\beta-1} ,
$$
where $\sigma$ is the slightly modified weight $\sigma ^2 (x) = 1+ x^2 /(2\beta-1)$. In other words, we have the optimal weighted Poincar\'e inequality: for all $f\in \CC ^\infty _0 (\rr)$,
$$
\frac{2(\beta-1)}{2\beta-1} \Var _{\mu_\beta} (f) \leq \int_\rr  \left( 1+ \frac{x^2}{2\beta-1} \right) f'(x) ^2 \mu_\beta (dx) .
$$
Passing through the limit as $\beta$ tends to infinity yields the classical Poincar\'e inequality for the standard Gaussian distribution, with the optimal constant 1.

\end{document}